\documentclass[11pt]{article}

\usepackage{amsmath,amssymb,latexsym,amsfonts,mathrsfs,graphicx,verbatim}

\usepackage{graphicx}
\usepackage{setspace}
\usepackage{program}
\usepackage{subfig}

\renewcommand{\quad}{$~~~\;\;\;$}
\newtheorem{theorem}{Theorem}
\newtheorem{remark}{Remark}
\newtheorem{proposition}{Proposition}
\newtheorem{lemma}{Lemma}
\newtheorem{corollary}{Corollary}
\newtheorem{example}{Example}

\newtheorem{conjecture}{Conjecture}

\newcommand{\conv}{\mbox{\rm conv}}

\newcommand{\transp}{{^{\rm T}}}

\newcommand{\vertiii}[1]{{\left\vert\kern-0.25ex\left\vert\kern-0.25ex\left\vert #1 
    \right\vert\kern-0.25ex\right\vert\kern-0.25ex\right\vert}}

\renewcommand{\R}{\mathbb{R}}

\newcommand{\matr}[1]{\begin{bmatrix} #1 \end{bmatrix}}    
\def\transp{^{\rm T}}

\newcommand{\ip}[2]{\left\langle #1 , #2 \right\rangle}    

\providecommand{\newoperator}[3]{%
  \newcommand*{#1}{\mathop{#2}#3}}

\newoperator{\argmax}{\mathrm{argmax}}{}
\newoperator{\argmin}{\mathrm{argmin}}{}

\newcommand{\dmin}{\displaystyle\min}
\newcommand{\dmax}{\displaystyle\max}
\newcommand{\dist}{\mathsf{dist}}

\usepackage{color}

\author{Javier Pe\~na\thanks{Tepper School of Business,
Carnegie Mellon University, USA, {\tt jfp@andrew.cmu.edu}}
\and Daniel Rodr\'iguez\thanks{Department of Mathematical Sciences, Carnegie
Mellon University, USA, {\tt drod@cmu.edu }}
\and Negar Soheili\thanks{College of Business Administration,  University of Illinois at Chicago, USA, {\tt nazad@uic.edu }}
}
\title{On the von Neumann and Frank-Wolfe Algorithms with Away Steps}

\begin{document}
\maketitle
\abstract{
The von Neumann algorithm is a simple coordinate-descent algorithm to determine whether the origin belongs to a polytope generated by a finite set of points. When the origin is in the {\em interior} of the polytope, the algorithm generates a sequence of points in the polytope that converges linearly to zero. The algorithm's rate of convergence depends on the radius of the largest ball around the origin contained in the polytope. 

\medskip

We show that under the weaker condition that the origin is in the polytope, possibly on its boundary, a variant of the von Neumann algorithm that includes {\em away steps} generates a sequence of points in the polytope that converges linearly to zero.
The new algorithm's rate of convergence depends on a certain geometric parameter of the polytope that extends the above radius but is always positive.  Our linear convergence result and geometric insights also extend to a variant of the Frank-Wolfe algorithm with away steps for minimizing a convex quadratic function over a  polytope.
}  

\newpage

\section{Introduction}

Assume $A = \matr{a_1 & \cdots & a_n} \in \R^{m\times n}$ with $\|a_i\|_2 =1, \; i=1,\dots,n.$  The von Neumann algorithm, communicated by von Neumann to Dantzig in the late 1940s and discussed later by Dantzig in an unpublished manuscript~\cite{Dant92}, is a simple algorithm to solve the feasibility problem:
\begin{center} Is $0 \in \conv(A) = \conv\{a_1,\dots,a_n\}$?  
\end{center}
More precisely, the algorithm aims to
find an approximate solution to the problem
\begin{equation}\label{von.neumann}
Ax = 0,  \; x\in \Delta_{n-1} = \{x\in \R^n_+: \|x\|_1=1\}.
\end{equation}
The algorithm starts from an arbitrary point $x_0\in \Delta_{n-1}$. At the $k$-th iteration the algorithm updates the current trial solution $x_k\in \Delta_{n-1}$ as follows.  First, it finds the column $a_j$ of $A$ that forms the widest angle with $y_k:=Ax_k$.  If this angle is acute, i.e., $A\transp y_k > 0$, then the algorithm halts as the vector $y_k$ separates the origin from $\conv(A)$.  Otherwise the algorithm chooses $x_{k+1}\in\Delta_{n-1}$ so that $y_{k+1}:=Ax_{k+1}$ is the minimum-norm convex combination of $Ax_k$ and $a_j$.  Let $e_j \in \Delta_{n-1}$ denote the $n$-dimensional vector with $j$-th component equal to one and all other components equal to zero.  To ease notation, we shall write $\|\cdot\|$ for $\|\cdot\|_2$ throughout the paper.

\bigskip
\noindent
{\bf Von Neumann Algorithm}

\begin{enumerate}
\item pick $x_0 \in \Delta_{n-1}$;  put $y_0:=Ax_0;\; k:=0$.
\item for $k=0,1,2,\dots$\\
\quad if $A\transp y_k > 0$ then HALT: $0\not \in\conv(A)$\\
\quad $j := \displaystyle\argmin_{i=1,\dots,n} \ip{a_i}{y_k};$\\
\quad $\theta_k := \displaystyle\argmin_{\theta \in [0,1]}\{\|y_k + \theta (a_j-y_k)\|\};$\\
\quad $x_{k+1} := (1-\theta_k)x_k + \theta_k e_j;\;\; y_{k+1} := Ax_{k+1};$ \\
end for 
\end{enumerate}
The von Neumann algorithm can be seen as a kind of coordinate-descent method for finding a solution to \eqref{von.neumann}:   At each iteration the algorithm judiciously selects a coordinate $j$ 
and {\em increases} the weight of the $j$-th component of $x_k$ while decreasing all of the others via a line-search step.  Like other currently popular coordinate-descent and first-order methods for convex optimization, the main attractive features of the von Neumann algorithm are its simplicity and low computational cost per iteration.  
Another attractive feature is its convergence rate.  Epelman and Freund~\cite{EpelF00} showed that the speed of convergence of the von Neumann algorithm can be characterized in terms of the following {\em condition measure} of the matrix $A$:
\begin{equation}\label{def.rho}
\rho(A):=\max_{z\in \R^m, \|z\|=1} \min_{i=1,\dots,n}\ip{a_i}{z}.
\end{equation}
The condition measure $\rho(A)$ was introduced by Goffin~\cite{Goff80} and later independently studied by Cheung and Cucker~\cite{CheuC01}.  The latter set of authors showed that  $\vert \rho(A)\vert$ is also a certain {\em distance to ill-posedness} in the spirit introduced and developed by Renegar~\cite{Rene95a,Rene95b}.  

Observe that $\rho(A)$ can also be written as
\begin{equation}\label{def.rho.alt}
\rho(A)=\max_{z\in \R^m, \|z\|=1} \min_{v\in\Delta_{n-1}} \ip{A\transp z}{v}=\max_{z\in \R^m, \|z\|=1} \min_{v\in\Delta_{n-1}} \ip{z}{Av}.
\end{equation}
Hence $\rho(A) > 0$ if and only if $0 \not \in \conv(A)$ and $\rho(A) < 0$ if and only if $0 \in \text{int}(\conv(A))$. When $\rho(A)>0$, this condition measure is closely related to the concept of margin in binary classification~\cite{Vapn98} and with the minimum enclosing ball problem in computational geometry~\cite{Clar10}.  The quantity $\rho(A)$ also has the following geometric interpretation  as discussed in \cite[Proposition 6.28]{BurgC13}. 
If $\rho(A) > 0$ then from~\eqref{def.rho.alt} and Lagrangian duality we get
\begin{equation}\label{rho.positive}
\begin{array}{rcl}
\rho(A) &=& \dmax_{z\in \R^m, \|z\|\le 1} \dmin_{v\in\Delta_{n-1}} \ip{z}{Av} \\&=& \dmin_{v\in\Delta_{n-1}} \dmax_{z\in \R^m, \|z\|\le 1}\ip{z}{Av} \\&=& \dmin\{\|y\|: y \in \conv(A)\} \\ &=& \dist(0,\partial\conv(A)).
\end{array}
\end{equation}
On the other hand, if $\rho(A) \le 0$ then \eqref{def.rho.alt} yields
\begin{equation}\label{rho.negative}
\begin{array}{rcl}
\vert  \rho(A) \vert &=& - \rho(A) \\ &=&  
\dmin_{z\in \R^m, \|z\|=1} \dmax_{v\in\Delta_{n-1}} \ip{z}{Av} \\
&=& 
\dmax\{r: \|y\| \le r \Rightarrow y \in \conv(A)\} \\&=& \dist(0,\partial\conv(A)).
\end{array}
\end{equation}
In either case $\vert\rho(A)\vert = \dist(0,\partial\conv(A)).$  Furthermore, observe that under the assumption $A = \matr{a_1 & \cdots & a_n} \in \R^{m\times n}$ with $\|a_i\| =1, \; i=1,\dots,n$ it follows that $\vert\ip{z}{Av}\vert\le 1$ for all $z\in \R^m, \|z\| =1$ and $v\in\Delta_{n-1}$.  In particular, from \eqref{def.rho.alt} it follows that $\vert\rho(A)\vert \le 1.$

\medskip

Epelman and Freund~\cite{EpelF00} showed the following properties of the von Neumann algorithm.  When $\rho(A) < 0$ the algorithm generates iterates $x_k\in \Delta_{n-1}, \; k=1, 2,\dots$ such that 
\begin{equation}\label{conv.von.neumann}
\|Ax_k\|^2 \le \left(1-\rho(A)^2\right)^k\|Ax_0\|^2.
\end{equation}
On the other hand, the iterates $x_k\in \Delta_{n-1}$ also satisfy
$
\|Ax_k\|^2 \le \frac{1}{k}
$
as long as the algorithm has not halted.  In particular, if $\rho(A) > 0$ then by \eqref{rho.positive} the algorithm must halt with a certificate of infeasibility $A\transp y_k >0 $ for $0\not\in \conv(A)$ in at most $\frac{1}{\rho(A)^2}$ iterations.  The latter bound is identical to a classical convergence bound for the perceptron algorithm~\cite{Bloc62,Novi62}.  This is not a coincidence as there is a nice duality between the perceptron and the von Neumann algorithms~\cite{LiT13,SoheP13}.

We show that a variant of the von Neumann algorithm with {\em away steps}  has the following stronger convergence properties.  When $0\in \conv(A)$, possibly on its boundary, the algorithm generates a sequence $x_k\in \Delta_{n-1}, \; k=1,2,\dots$ satisfying
\begin{equation}\label{conv.von.neumann.away}
\|Ax_k\|^2 \le \left(1-\frac{w(A)^2}{16}\right)^{k/2}\|Ax_0\|^2.
\end{equation}
The quantity $w(A)$ is a kind of {\em relative width} of $\text{conv}(A)$ that is at least as large as $\vert\rho(A)\vert$.  However, unlike $\vert\rho(A)\vert$ the relative width $w(A)$ is positive for any non-zero matrix $A\in \R^{m\times n}$ provided $0\in \conv(A)$.  When $\rho(A) >0$, or equivalently $0\not \in\conv(A)$, the von Neumann algorithm with away steps finds a certificate of infeasibility $A\transp y_k > 0$ for $0\not\in \conv(A)$ in at most $\frac{8}{\rho(A)^2}$ iterations.

\begin{figure}[h!]
\centering
{\includegraphics[width=0.725\textwidth]{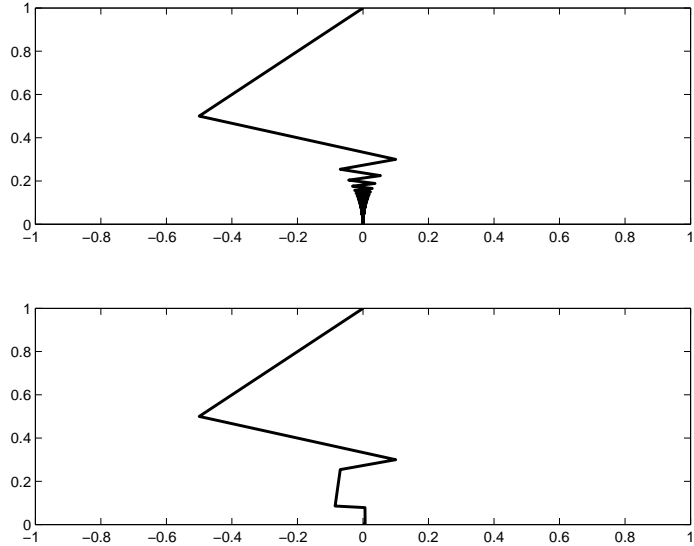}}
\caption{Iterates generated by the von Neumann algorithm (top) and its variant with away steps (bottom)}
\label{the.fig}
\end{figure}

Figure~\ref{the.fig} illustrates the different behavior of the von Neumann algorithm and the variant with away steps described in Section~\ref{sec.von.neumann.away} for $A = \matr{1 & 0 & 0 \\ 0 & -1 & 1}$.  The figure depicts the path of iterates $\{y_k:k=0,1,\dots\}$ generated by each algorithm starting from $y_0 = \matr{1\\0}.$ The zig-zagging behavior in the first case occurs because after the third iteration the search direction is nearly perpendicular to the current iterate and as a consequence the algorithm makes slow progress.  By contrast, in the second case the away steps provide alternative search directions that enable the algorithm to make faster progress.

The von Neumann algorithm can be seen as a special case of the Frank-Wolfe (also known 
as conditional gradient) algorithm~\cite{FranW56,Jagg13}.  The von Neumann algorithm is also nearly identical to an algorithm for minimizing a quadratic form over a convex set independently developed by Gilbert~\cite{Gilb66}.  The name ``Gilbert's algorithm''  appears to be more popular in the computational geometry literature~\cite{GartJ09}. 

We show that a linear convergence result similar to~\eqref{conv.von.neumann.away} also holds for a version of the Frank-Wolfe algorithm  with away steps for minimizing a strongly convex quadratic function over a polytope.  This variant of the Frank-Wolfe algorithm with away steps was introduced by Wolfe~\cite{Wolf70} and has been subsequently studied by various authors.  In particular,  linear convergence results similar to ours have been previously established in~\cite{GarbH13,GuelM86,Jagg13,LacoJ15} and more recently in~\cite{BeckS15}.  Linear convergence results in the same spirit also hold for the randomized Kaczmarz algorithm~\cite{StroV09} and for the methods of randomized coordinate descent and iterated projections~\cite{LeveL10}.  The computational article~\cite{GoncSG09} also reports numerical experiments for variants of the von Neumann algorithm with away steps. Our main contributions are the succinct and transparent proofs of linear convergence results that highlight the role of the relative width $w(A)$ and a closely related {\em restricted width} $\phi(A)$.  Our presentation unveils a deep connection between problem conditioning as encompassed by the quantities $w(A), \phi(A)$ and the behavior of the von Neumann and Frank-Wolfe algorithms with away steps.  We also provide some lower bounds on $w(A)$ and $\phi(A)$ in terms of certain radii quantities  that naturally extend $\rho(A)$.  We note that the linear convergence results in \cite{LacoJ15} are stated in terms of a certain {\em pyramidal width} whose geometric intuition and properties appear to be less understood than those of $w(A)$ and $\phi(A)$.  

The rest of the paper is organized as follows.  In Section~\ref{sec.von.neumann.away} we describe a von Neumann Algorithm with Away Steps and establish its main convergence result in terms of the relative width $w(A)$.  Section~\ref{sec.quadratic} extends our main result to the more general problem of minimizing a quadratic function over the polytope $\text{conv}(A)$.  Finally, Section~\ref{sec.width} discusses some properties of the relative and restricted widths.

\section{Von Neumann Algorithm with Away Steps}\label{sec.von.neumann.away}

Throughout this section we assume $A = \matr{a_1 & \cdots & a_n} \in \R^{m\times n}$ with $\|a_i\| = 1, \; i=1,\dots,n.$   We next consider a variant of the von Neumann Algorithm that includes so-called ``away'' steps.  To that end, at each iteration, in addition to a ``regular step'' the algorithm considers an alternative ``away step''.  Each of these away steps  identifies a coordinate $\ell$ such that the $\ell$-th component of $x_k$ is positive and {\em decreases} the weight of the $\ell$-th component of $x_k$.
The algorithm needs to keep track of the {\em support}, that is, the set of positive entries of a vector.  To that end,  given $x\in \R^n_+$, let the support of $x$ be defined as $$S(x):=\{i\in \{1,\dots,n\}: x_i > 0\}.$$

\bigskip

\noindent
{\bf Von Neumann Algorithm with Away Steps}

\begin{enumerate}
\item pick $x_0 \in \Delta_{n-1}$; put $y_0:= Ax_0;\; k:=0;$.
\item for $k=0,1,2,\dots$\\
\quad if $A\transp y_k > 0$ then HALT:  $0\not \in\conv(A)$\\
\quad $j := \displaystyle\argmin_{i=1,\dots,n} \ip{a_i}{y_k};\;\;\ell := \displaystyle\argmax_{i\in S(x_k)} \ip{a_i}{y_k};$\\
\quad if $\ip{a_j - y_k}{y_k} < \ip{y_k-a_\ell}{y_k}$
 then  (regular step)\\
\quad \quad  $a:=a_j - y_k;\; u:=e_j - x_k;\;\theta_{\max} := 1$ \\
\quad else (away step)\\
\quad \quad  $a:= y_k - a_\ell; \; u:= x_k -e_{\ell};\;\theta_{\max} := \frac{(x_k)_\ell}{1-(x_k)_\ell}$ \\
\quad endif \\
\quad $\theta_k := \displaystyle\argmin_{\theta \in [0,\theta_{\max}]}\{\|y_k + \theta a\|\};$\\
\quad $ x_{k+1} := x_k + \theta_k u; \;\; y_{k+1} := Ax_{k+1}$\\
end for 
\end{enumerate}

Note that the above von Neumann Algorithm with Away Steps can also be applied to any non-zero matrix $A = \matr{a_1 & \cdots & a_n}$. The assumption that  the columns of $A$ are normalized, i.e., $\|a_i\|=1, \; i=1,\dots,n,$ simplifies our notation and exposition.  In Section~\ref{sec.quadratic} below we extend our discussion to the case when the columns of $A$ are not necessarily normalized.

Observe that the iterates $x_k, k=0,1,\dots,$ generated by the above von Neumann Algorithm with Away Steps satisfy $x_k \in \Delta_{n-1}$.  This fact follows by induction:  By construction, $x_0 \in \Delta_{n-1}$.  At iteration $k$ we have $x_{k+1} = x_k + \theta_k u$ where $x_k \in \Delta_{n-1}$ and the components of $u$ add up to zero as $u$ is either $e_j-x_k$ or $e_{\ell} - x_k$.  The bound $\theta_k \in[0,\theta_{\max}] $ in turn guarantees that $x_{k+1} \ge 0$ and so $\|x_{k+1}\|_1 = \|x_k\|_1 = 1.$

\bigskip

Define the relative width $w(A)$  of conv$(A)$ as 
\begin{equation}\label{def.rel.width}
w(A):=\dmin_{x\ge 0, Ax\ne 0} \dmax_{\ell,j}\left\{\frac{\ip{Ax}{a_\ell-a_j}}{\|Ax\|}: \ell\in S(x), \; j\in\{1,\dots,n\}
\right\}.
\end{equation}
The next proposition shows that $w(A) \ge \vert\rho(A)\vert$ when $0\in \conv(A)$.  To that end, observe that  $w(A)$  can also be written as
\begin{equation}\label{def.rel.width.alt}
w(A) = \dmin_{x\ge 0, Ax\ne 0} \dmax_{u,v}\left\{\frac{\ip{Ax}{Au-Av}}{\|Ax\|}: u,v\in\Delta_{n-1}, S(u) \subseteq S(x) \right\}.
\end{equation}

\begin{proposition} If $A$ is such that $0\in \conv(A)$ then $w(A) \ge \vert\rho(A)\vert$.
\end{proposition}
\proof
Since $0\in \conv(A)$, equation \eqref{def.rho.alt} yields
\[
\rho(A) = \max_{z\in \R^m, \|z\|=1} \min_{v\in\Delta_{n-1}} \ip{z}{Av} \le 0.
\]
In particular,
\begin{equation}\label{eq.abs.rho}
\begin{array}{rcl}
\vert\rho(A)\vert &=& \dmin_{z\in \R^m, \|z\|=1} \max_{v\in\Delta_{n-1}} \ip{z}{-Av} \\
 & \le & \dmin_{x\ge 0, Ax\ne 0} \dmax_{v\in \Delta_{n-1}}\frac{\ip{Ax}{-Av}}{\|Ax\|}.
 \end{array}
\end{equation}
Hence from \eqref{def.rel.width.alt} we get
\[ 
w(A) \ge  \dmin_{x\ge 0, Ax\ne 0} \dmax_{v\in \Delta_{n-1}}\frac{\ip{Ax}{-Av}}{\|Ax\|}  
 \ge  \vert\rho(A)\vert.
\]
The first inequality holds because we can choose $u = \frac{x}{\|x\|_1}$ in \eqref{def.rel.width.alt}. The second inequality follows from \eqref{eq.abs.rho}.
\qed

Observe that under the assumption $A = \matr{a_1 & \cdots & a_n} \in \R^{m\times n}$ with $\|a_i\| =1, \; i=1,\dots,n$ it follows that $\|Au - Av\|\le 2$ for all $u,v\in\Delta_{n-1}$.  In particular, from \eqref{def.rel.width.alt} it follows that $w(A) \le 2.$  In Section~\ref{sec.width} below we discuss some additional properties of $w(A)$.  In particular, we will formally prove that $w(A)>0$ for any nonzero matrix $A\in \R^{m\times n}$ such that $0 \in \conv(A)$.

We are now ready to state the main properties of the von Neumann algorithm with away steps.

\begin{theorem}~\label{lin.convergence.vn} Assume $x_0 \in \Delta_{n-1}$ is one of the extreme points of $\Delta_{n-1}$.
\begin{description}
\item[(a)] If  $0\in \text{conv}(A)$ then the iterates $x_k\in \Delta_{n-1}, y_k=Ax_k, \; k=1,2,\dots$ generated by the von Neumann Algorithm with Away Steps satisfy
\[
\|y_k\|^2 \le \left(1-\frac{w(A)^2}{16}\right)^{k/2} \|y_0\|^2. 
\]
\item[(b)] The iterates $x_k\in \Delta_{n-1}, y_k=Ax_k, \; k=1,2,\dots$ generated by the von Neumann Algorithm with Away Steps also satisfy
\[
\|y_k\|^2 \le \frac{8}{k}
\]
as long as the algorithm has not halted. In particular, if $0\not\in \text{conv}(A)$ then the von Neumann Algorithm with Away Steps finds a certificate of infeasibility $A\transp y_k > 0$ for $0\not\in \conv(A)$ in at most $\frac{8}{\rho(A)^2}$ iterations.
\end{description}
\end{theorem}
The crux of the proof of Theorem~\ref{lin.convergence.vn} is the following elementary lemma.

\begin{lemma}\label{the.lemma.vn}
Assume $a,y\in\R^m$ satisfy $\ip{a}{y} <0$.  Then 
\[
\min_{\theta \ge 0} \|y + \theta a\|^2 = \|y\|^2 - \frac{\ip{a}{y}^2}{\|a\|^2},
\]
and the minimum is attained at $\theta = -\frac{\ip{a}{y}}{\|a\|^2}.$
\end{lemma}

\bigskip

\noindent
{\bf Proof of Theorem~\ref{lin.convergence.vn}:}

\begin{description}
\item[(a)] The algorithm generates $y_{k+1}$ by solving a problem of the form
\[
\|y_{k+1}\|^2 = \min_{\theta\in[0,\theta_{\max}]}\|y_k + \theta a\|^2
\]
where $a=a_j-y_k$ or $a=y_k - a_{\ell}$ is chosen so that $\ip{a}{y_k} = \min\{\ip{y_k-a_{\ell}}{y_k},\ip{a_j-y_k}{y_k}\}.$  In particular,
\begin{equation}\label{eq.a}
\begin{array}{rcl}
-\ip{a}{y_k}
& \ge&  \displaystyle\frac{1}{2}\left(\ip{a_{\ell} - y_k}{y_k}+\ip{y_k- a_j}{y_k}\right) \\[2ex]
&=& \displaystyle\frac{1}{2}\ip{a_{\ell} - a_j}{y_k} \\[2ex]
&\ge& \displaystyle\frac{1}{2} w(A) \|y_k\|.\end{array}
 \end{equation}
If $\theta_k < \theta_{\max}$ then Lemma~\ref{the.lemma.vn} applied to $y:=y_k$ yields
\[
\|y_{k+1}\|^2 = \|y_k\|^2 - \frac{\ip{a}{y_k}^2}{\|a\|^2} \le \|y_k\|^2 - \frac{w(A)^2}{16} \|y_k\|^2.
\]
The second inequality follows from \eqref{eq.a} and $\|a\| \le 1 + \|y_k\| \le 2$.
Thus each time the algorithm performs an iteration with $\theta_k<\theta_{\max}$, the value of $\|y_k\|^2$ decreases at least by the factor $1-\frac{w(A)^2}{16}$.  To conclude, it suffices to show that after $N$ iterations the number of iterations with $\theta_k < \theta_{\max}$ is at least $N/2$.  To that end, we apply the following argument from \cite{LacoJ15}:  
Observe that when $\theta_k = \theta_{\max}$ we have $\vert S(x_{k+1})\vert<\vert S(x_k)\vert$.   On the other hand, when $\theta_k < \theta_{\max}$ we have $\vert S(x_{k+1})\vert\le \vert S(x_k)\vert+1$.  
Since $\vert S(x_0)\vert=1$ and $\vert S(x)\vert\ge 1$ for every $x\in \Delta_{n-1}$, after any number of iterations there must have been at least as many iterations with $\theta_k < \theta_{\max}$ as there have been iterations with $\theta_k = \theta_{\max}.$ Hence after $N$ iterations, the number of iterations with $\theta_k < \theta_{\max}$ is at least $N/2$.
\item[(b)] Proceed as above but note that if the algorithm does not halt at the $k$-th iteration  then $\ip{a}{y_k} \le \ip{a_j-y_k}{y_k} \le -\|y_k\|^2$.  Thus each time the algorithm performs an iteration with $\theta_k<\theta_{\max}$, we have
\begin{equation}\label{eq.induction}
\|y_{k+1}\|^2 \le \|y_k\|^2 - \frac{\ip{a}{y_k}^2}{\|a\|^2}\le\|y_k\|^2 - \frac{\|y_k\|^4}{4}.
\end{equation}
Assume the algorithm has not halted after $N$ iterations.  Let $m$ be the number of iterations with $\theta_k < \theta_{\max}$ up to iteration $N$.  If $\|y_N\|^2 \le \frac{4}{m}$ and $\theta_{N} < \theta_{\max}$ then from \eqref{eq.induction} we get
\[
\|y_{N+1}\|^2 \le \frac{4}{m} - \frac{4}{m^2} = \frac{4(m-1)}{m^2} \le \frac{4}{m+1}.
\]
It follows by induction that if the algorithm has not halted after $N$ iterations then $\|y_N\|^2 \le \frac{4}{m}$.  As in part (a), it must be the case that $m \ge \frac{N}{2}$ and consequently $\|y_N\|^2\le \frac{8}{N}.$ Finally, if $0\not\in \text{conv}(A)$ then $\rho(A) = \min\{\|y\|: y\in\conv(A)\}>0$ and so the algorithm must halt with a certificate of infeasibility $A\transp y_k > 0$ 
for $0\not\in \text{conv}(A)$ after at most $\frac{8}{\rho(A)^2}$ iterations.
\end{description}
\qed

\section{Frank-Wolfe Algorithm with Away Steps}\label{sec.quadratic}

Throughout this section assume $A = \matr{a_1 & \cdots & a_n} \in \R^{m\times n}$ is a non-zero matrix, and $f(y) = \displaystyle\frac{1}{2}\ip{y}{Qy} + \ip{b}{y}$ for a symmetric positive definite matrix $Q\in\R^{m\times m}$ and $b \in \R^m$.  Consider the problem
\begin{equation}\label{min.quad}
\min_{y\in \text{conv}(A)} f(y) \Leftrightarrow 
\min_{x\in \Delta_{n-1}} f(Ax). \end{equation}
Observe that in contrast to Section \ref{sec.von.neumann.away}, we do not assume that the columns of $A$ are normalized in this section. 

Problem \eqref{von.neumann} can be seen as a special case of \eqref{min.quad} when $Q = I$ and $b = 0$.  The von Neumann Algorithm can also be seen as a special case of the Frank-Wolfe Algorithm~\cite{FranW56} for \eqref{min.quad}.  This section extends the ideas and results from Section~\ref{sec.von.neumann.away} to the following variant of the Frank-Wolfe algorithm with away steps.  This variant can be traced back to Wolfe~\cite{Wolf70}. It has been a subject of study in a number of papers~\cite{BeckS15,GarbH13,GoncSG09,GuelM86,LacoJ15}.

\bigskip

\noindent
{\bf Frank-Wolfe Algorithm with Away Steps}

\begin{enumerate}
\item pick $x_0 \in \Delta_{n-1}$; put $y_0:= Ax_0;\; k:=0;$.
\item for $k=0,1,2,\dots$\\
\quad $j := \displaystyle\argmin_{i=1,\dots,n} \ip{a_i}{\nabla f(y_k)};\; \ell := \displaystyle\argmax_{i\in S(x_k)}  \ip{a_i}{\nabla f(y_k)};$\\
\quad if $\ip{a_j-y_k}{\nabla f(y_k)} < \ip{y_k-a_{\ell}}{\nabla f(y_k)}$ 
 then  (regular step)\\
\quad \quad  $a:=a_j - y_k;\; u:=e_j - x_k; \; \theta_{\max} := 1$ \\
\quad else (away step)\\
\quad \quad  $a:= y_k - a_\ell; \; u:= x_k -e_{\ell}; \; \theta_{\max} := \frac{(x_k)_\ell}{1-(x_k)_\ell}$ \\
\quad endif \\
\quad $\theta_k := \displaystyle\argmin_{\theta \in [0,\theta_{\max}]}f(y_k + \theta a)$\\
\quad $x_{k+1} := x_k + \theta_k u; \; \; y_{k+1} := Ax_{k+1}$\\
end for 
\end{enumerate}
Observe that the computation of $\theta_k$ in the second to last step reduces to minimizing a one-dimensional convex quadratic function over the interval $[0,\theta_{\max}]$.

\medskip

We next present a general version of Theorem~\ref{lin.convergence.vn} for the above Frank-Wolfe Algorithm with Away Steps.  The linear convergence result depends on a certain restricted width and diameter defined as follows.  
For $x\ge 0$ with $Ax \ne 0$ let
\begin{multline*}
\phi(A,x):=\\ \sup\left\{\lambda > 0: \exists u,v\in \Delta_{n-1}, \; S(u) \subseteq S(x),\; Au - Av = \frac{\lambda}{\|Ax\|}Ax\right\}.
\end{multline*}
Define the restricted width $\phi(A)$ and diameter $d(A)$ of conv$(A)$ as follows.
\begin{equation}\label{def.width}
\phi(A):= \dmin_{x}\left\{ \phi(A,x):  x\ge 0,\;Ax \ne 0\right\},
\end{equation}
and
\begin{equation}\label{diameter}
d(A):= \max_{x,u\in \Delta_{n-1}}\|Ax - Au\|.
\end{equation}
Observe that for $x\ge 0$ with $Ax \ne 0$ 
\[
\phi(A,x)
\le \dmax_{u,v}\left\{\frac{\ip{Ax}{Au-Av}}{\|Ax\|}: u,v\in\Delta_{n-1}, S(u) \subseteq S(x) \right\}.
\]
Thus \eqref{def.rel.width.alt} and \eqref{def.width} imply that $w(A) \ge \phi(A)$ for all nonzero $A\in \R^{m\times n}$. Furthermore, the restricted width $\phi(A)$ can be seen as an extension of the radius $\rho(A)$ defined in \eqref{def.rho}. Indeed, when $0\in \text{int}(\conv(A))$, we have $\text{span}(A) = \R^m$.  Hence \eqref{rho.negative}  can alternatively be written as
\[
\vert\rho(A)\vert:=\min_{x\ge 0, Ax \ne 0} \max\left\{\lambda: \exists v \in \Delta_{n-1}, \, -Av = \frac{\lambda}{\|Ax\|}Ax \right\}.
\]
This implies that $\phi(A,x) \ge \vert\rho(A)\vert + \frac{\|Ax\|}{\|x\|_1}$ for all $x\ge 0$ with $Ax \ne 0$.  Hence the following inequality readily follows $$\phi(A) \ge \vert\rho(A)\vert.$$

Section~\ref{sec.width} presents a stronger lower bound on $\phi(A)$ in terms of certain variants of $\rho(A)$.  In particular, we will show that  $\phi(A)>0$, and consequently $w(A)>0$, for any nonzero matrix $A\in \R^{m\times n}$ such that $0 \in \conv(A)$.

The linear convergence property of the von Neumann algorithm with away steps, as stated in Theorem~\ref{lin.convergence.vn}(a),  extends as follows.

\begin{theorem}~\label{lin.convergence.quad}  Assume $x^*\in\Delta_{n-1}$ is a minimizer of \eqref{min.quad}.  Let $y^*=Ax^*$ and $\bar A := Q^{1/2}\matr{a_1- y^* & \cdots & a_n-y^*}.$ If $x_0 \in \Delta_{n-1}$ is one of the extreme points of $\Delta_{n-1}$ then the iterates $x_k\in \Delta_{n-1}, y_k = Ax_k, \; k=1,2,\dots$ generated by the Frank-Wolfe Algorithm with Away Steps satisfy
\begin{equation}\label{norm.bound.gral}
f(y_k) - f(y^*) \le \left(1-\frac{\phi(\bar A)^2}{4d(\bar A)^2}\right)^{k/2} (f(y_0) - f(y^*)).
\end{equation}
\end{theorem}

The proof of Theorem~\ref{lin.convergence.quad} relies on the following two lemmas.  The first one is similar to Lemma~\ref{the.lemma.vn} and also follows via a straightforward calculation.

\begin{lemma}\label{the.lemma}
Assume $f$ is as above and $a,y\in\R^m$ satisfy $\ip{a}{\nabla f(y)} <0$.  Then 
\[
\min_{\theta \ge 0} f(y + \theta a) = f(y) - \frac{\ip{a}{\nabla f(y)}^2}{2\ip{a}{Qa}},
\]
and the minimum is attained at $\theta = -\frac{\ip{a}{\nabla f(y)}}{\ip{a}{Qa}}.$
\end{lemma}

\begin{lemma}\label{opt.cond}
Assume $f, A, y^*, \bar A$ are as in Theorem~\ref{lin.convergence.quad} above.  Then for all $x \in \Delta_{n-1}$ 
\[
\max_{\ell\in S(x), j=1,\dots,n}\ip{\nabla f(Ax)}{a_\ell - a_j} \ge \phi(\bar A)\sqrt{2(f(Ax) - f(y^*))}.
\]
\end{lemma}
\proof
Let $y := Ax \in \conv(A)$.  Assume $y\ne y^*$ as otherwise there is nothing to show. Since $y^*$ minimizes \eqref{min.quad},  we have $\ip{\nabla f(y^*)}{y-y^*} \ge 0$.  For ease of notation put $\delta:=\ip{\nabla f(y^*)}{y-y^*}$ and 
$\|y-y^*\|_Q^2:= \ip{y-y^*}{Q(y-y^*)}.$ It readily follows that
\[
\ip{\nabla f(y)}{y-y^*}  = \|y-y^*\|_Q^2 + \delta \ge 0,
\]
and
\[
2(f(y) - f(y^*)) = \|y-y^*\|_Q^2 + 2\delta.
\]
Hence,
\begin{align*}
\ip{\nabla f(y)}{y-y^*}^2  &= (\|y-y^*\|_Q^2 + \delta)^2 \\&\ge \|y-y^*\|_Q^2(\|y-y^*\|_Q^2+2\delta) \\&= 2\|y-y^*\|_Q^2(f(y) - f(y^*)).
\end{align*}
Thus
\begin{equation}\label{eqn.1}
\frac{\ip{\nabla f(y)}{y-y^*}}{\|y-y^*\|_Q} \ge \sqrt{2(f(y) - f(y^*))}. 
\end{equation}
On the other hand, by the definition of $\phi(A)$ there exist $u,v\in\Delta_{n-1}$ with $S(u) \subseteq S(x)$
and $\lambda \ge \phi(\bar A)$ such that $
\bar Au - \bar A v = \frac{\lambda}{\|\bar Ax\|} \bar A x .$  Since $\bar A x = Q^{1/2}(Ax-y^*) = Q^{1/2}(y-y^*)$,
the latter equation can be rewritten as
\begin{equation}\label{eqn.2}
Au - Av = \frac{\lambda}{\|y-y^*\|_Q} (y-y^*).
\end{equation}
Putting \eqref{eqn.1} and \eqref{eqn.2} together we get\[
\ip{\nabla f(y)}{Au - Av} =  \frac{\lambda \ip{\nabla f(y)}{y - y^*}}{\|y-y^*\|_Q} \\ 
\ge \phi(\bar A) \sqrt{2(f(y) - f(y^*))}.
\]
To finish, observe that
\begin{align*}
\max_{\ell\in S(x), j=1,\dots,n}\ip{\nabla f(Ax)}{a_\ell - a_j} &\ge 
\ip{\nabla f(y)}{Au - Av} \\ 
 &\ge \phi(\bar A) \sqrt{2(f(Ax) - f(y^*))}.
\end{align*}
\qed

\noindent{\bf Proof of Theorem~\ref{lin.convergence.quad}:}  This is a modification of the proof of Theorem~\ref{lin.convergence.vn}(a).  At iteration $k$ the algorithm yields $y_{k+1}$ such that
\[
f(y_{k+1}) = \min_{\theta\in[0,\theta_{\max}]}f(y_k + \theta a)
\]
where $a=a_j-y_k$ or $a=y_k - a_{\ell}$, and 
$$-\ip{\nabla f(y_k)}{a} > \frac{1}{2}\ip{\nabla f(y_k)}{a_{\ell} - a_j} \ge \frac{1}{2} \phi(\bar A) \sqrt{2(f(y_k)-f(y^*)}.$$
The second inequality above follows from Lemma~\ref{opt.cond}.
If $\theta_k < \theta_{\max}$ then Lemma~\ref{the.lemma} applied to $y:=y_k$ yields
\[
f(y_{k+1}) = f(y_k) - \frac{\ip{a}{\nabla f(y_k)}^2}{2\ip{a}{Qa}} \le 
f(y_k) - \frac{\phi(\bar A)^2}{4d(\bar A)^2}(f(y_k)-f(y^*)).
\]
That is,
\[
f(y_{k+1}) - f(y^*) \le \left(1- \frac{\phi(\bar A)^2}{4d(\bar A)^2}\right)(f(y_k)-f(y^*)).\]
Then proceeding as in the last part of the proof of Theorem~\ref{lin.convergence.vn}(a) we obtain~\eqref{norm.bound.gral}.
\qed

\begin{remark}\label{the.remark}
A closer look at the proof of Theorem~\ref{lin.convergence.quad} reveals that the convergence bound~\eqref{norm.bound.gral} can be sharpened as follows: Replace $\phi(\bar A)$ with $w_f(A)\ge \phi(\bar A)$, where $w_f(A)$ is the following extension of $w(A)$:
\begin{multline*}
w_f(A):= \\
\dmin_{\tiny\begin{array}{c}x\in \Delta_{n-1}\\ Ax\ne y^*\end{array}} \dmax_{\ell,j}\left\{ 
\frac{\ip{\nabla f(Ax)}{a_\ell - a_j}}{\sqrt{2(f(Ax) - f(y^*))}}:
 \ell\in S(x), j\in\{1,\dots,n\}\right\}.
\end{multline*}
In the special case when $Q= I, b = 0$ problem~\eqref{min.quad} specializes to problem \eqref{von.neumann}.  In this case if $0\in\conv(A)$ then we have $y^* = 0$ and $w_f(A) = w(A)$. Hence the sharpened version of Theorem~\ref{lin.convergence.quad} yields
\[
\|y_k\|^2 \le \left(1-\frac{w(A)^2}{4d(A)^2}\right)^{k/2} \|y_0\|^2. 
\]
If in addition the columns of $A$ are normalized then $d(A) \le 2$ and we recover the bound in Theorem~\ref{lin.convergence.vn}(a).
\end{remark}

We have the following related conjecture concerning $w(A)$ and $\phi(A)$.

\begin{conjecture}  If $A\in \R^{m\times n}$ is non-zero and $0\in\conv(A)$ then $\phi(A) = w(A)$.
\end{conjecture}

\section{Some properties of the restricted width}\label{sec.width}

Throughout this section assume $A\in \R^{m\times n}$ is a nonzero matrix. As we noted in Section~\ref{sec.quadratic} above, $w(A) \ge \phi(A)$ and $\phi(A) \ge \vert\rho(A)\vert$ when $0\in \text{int}(\conv(A))$.    Our next result establishes a stronger lower bound on $\phi(A)$ in terms of some quantities that generalize $\rho(A)$ to the case when $0\in \partial \text{conv}(A)$.  To that end, we recall some terminology and results from \cite{CheuCP09}.  Assume $A = \matr{a_1 & \cdots & a_n} \in \R^{m\times n}$ is a non-zero matrix.   Then there exists a unique partition $B\cup N = \{1,\dots,n\}$ such that both 
$A_B x_B = 0, \; x_B > 0$ and $A_N\transp y > 0, \; A_B\transp y = 0$ are feasible.  In particular, $B\ne \emptyset$ if and only if $0 \in \text{conv}(A)$.  Also $N \ne \emptyset$ if and only if $ 0 \not\in \text{relint}(\conv(A)).$  Furthermore, if $a_i = 0$ then $i \in B$.

The above canonical partition $(B,N)$ allows us to refine the quantity  $\rho(A)$ defined by \eqref{def.rho} as follows.  Let $L := \text{span}(A_B)$ and $L^\perp := \{v\in \R^m : \ip{v}{y} = 0 \; \text{ for all } \; y \in L\}$.  By convention, $L = \{0\}$ and $L^\perp = \R^m$ when $B = \emptyset.$
If $L\ne \{0\}$, let $\rho_B(A)$ be defined as
\[
\rho_B(A):= \max_{z\in L,\|z\| = 1} \min_{i\in B} \ip{a_i}{z}.\]
Observe that if $B \ne \emptyset$, then $L = \{0\}$ only when $a_i = 0$ for all $i\in B$.

If $N\ne \emptyset$, let $\rho_N(A)$ be defined as
\[
\rho_N(A):= \max_{z\in L^\perp,\|z\| = 1} \min_{i\in N} \ip{a_i}{z}. \]
When $L\ne\{0\}$, it can be shown~\cite{CheuCP09} that 
 $\rho_B(A) < 0.$  Likewise, when $N\ne\emptyset$ it can be shown that $\rho_N(A)>0$.  In particular, the latter implies that 
\begin{equation}\label{rhoN}
\rho_N(A):= \max_{z\in L^\perp,\|z\| = 1} \min_{i\in N} \ip{a_i}{z} = 
\max_{z\in L^\perp,\|z\| \le 1} \min_{i\in N} \ip{a^\perp_i}{z},
\end{equation}
where $a^\perp_i$ is the orthogonal projection of $a_i$ onto $L^\perp$.
Let $A^\perp_N$ denote the matrix obtained by projecting each of the columns of $A_N$ onto $L^\perp.$  From \eqref{rhoN} and Lagrangian duality it follows that
\begin{equation}\label{rhoNdual}
\rho_N(A)=\dmin \{\|y\|: y \in \conv(A^\perp_N)\}.
\end{equation}
Similarly, it can be shown that if $L\ne\{0\}$ then
\begin{equation}\label{rhoBprimal}
\vert\rho_B(A)\vert = \max\{r: y\in L, \|y \|\le r \Rightarrow y\in \text{conv}(A_B)\}.
\end{equation}
Observe that \eqref{rhoNdual} and \eqref{rhoBprimal} nicely extend \eqref{rho.positive} and \eqref{rho.negative}.  Indeed, \eqref{rhoNdual} is identical to \eqref{rho.positive} when $B = \emptyset$.  Likewise, \eqref{rhoBprimal} is identical to \eqref{rho.negative} when $N = \emptyset$ and $\text{rank}(A) = m$.  Furthermore, \eqref{rhoNdual} and \eqref{rhoBprimal} imply that $\rho_N(A) = \dist(0,\partial \text{conv}(A^\perp_N))$ and $\vert\rho_B(A)\vert = \dist_L(0,\partial \text{conv}(A_B))$ thereby extending the fact that $\vert \rho(A)\vert = \dist(0,\partial\conv(A)).$

The next results show that $\phi(A)$ can be bounded below in terms of $\rho_B(A)$ and $\rho_N(A)$.  In particular, Corollary~\ref{bound.varrho} shows that $w(A)\ge \phi(A)>0$ whenever $A\ne 0$ and $0 \in \conv(A)$.

\begin{theorem}~\label{main.thm} 
Assume $A = \matr{a_1 & \cdots & a_n} \in \R^{m\times n}$ is a nonzero matrix.
\begin{description}
\item[(a)] 
If $N = \emptyset$ then $L \ne \{0\}$ and $\phi(A) \ge \vert\rho_B(A)\vert$.
\item[(b)] If $B= \emptyset$ then 
$
\phi(\bar A) \ge \rho_N(A)$ for $\bar A := \matr{A& 0}.$  
\item[(c)] If $B\ne \emptyset$ and $L = \{0\}$ then $\phi(A) \ge \rho_N(A)$.
\item[(d)]
If $N \ne \emptyset$ and $L \ne \{0\}$ then 
$
\phi(A) \ge \dfrac{\vert\rho_B(A)\vert\rho_N(A)}{\sqrt{\|A\|^2 + \rho_N(A)^2}}
,$ where $\|A\| = \dmax_{i=1,\dots,n} \|a_i\|.$
\end{description}
\end{theorem}
\proof
\begin{description}
\item[(a)] Assume $ x\ge 0$ is such that $y:= Ax \ne 0.$  In this case $y \in \text{span}(A_B) = L.$  Hence $L\ne\{0\}$ and by~\eqref{rhoBprimal} there exists $v \in \Delta_{n-1}$ and $r \ge \vert\rho_B(A)\vert$ such that $-A v =  \frac{r}{\|Ax\|}Ax$.  Thus for $u:=\frac{x}{\|x\|_1}$ we have $u,v\in \Delta_{n-1}$, $S(u)\subseteq S(x)$ and $Au - Av = \left(r + \frac{\|Ax\|}{\|x\|_1}\right)\frac{1}{\|Ax\|}Ax.$ It follows that $\phi(A,x) \ge r + \frac{\|Ax\|}{\|x\|_1} > \vert\rho_B(A)\vert$.  
\item[(b)] Assume $\bar x:=\matr{x\\t} \ge 0$ is such that $y:= \bar A \bar x = Ax \ne 0.$  
From~\eqref{rhoNdual} it follows that $\frac{\|Ax\|}{\|x\|_1} \ge \rho_N(A).$  Thus for 
$u := \matr{\frac{x}{\|x\|_1} \\ 0}, \; v:=e_{n+1}$ we have $u,v\in \Delta_{n-1}$, $S(u)\subseteq S(\bar x)$ and $\bar A u - \bar A v = \frac{\|Ax\|}{\|x\|_1}\frac{1}{\|Ax\|}Ax.$  It follows that $\phi(\bar A, \bar x)  \ge \frac{\|Ax\|}{\|x\|_1}\ge \rho_N(A).$ 

\item[(c)] Since $B \ne \emptyset$ and $L = \{0\}$, it follows that $A_B= 0$ and the columns of $A_N$ are precisely the non-zero columns of $A$.  Thus from part (b) we get $\phi\left(\matr{A_N & 0}\right) \ge \rho_N(A)$.  To finish, observe that 
$\phi(A) = \phi(\matr{A_N & 0})$ because $A_B = 0$.

\item[(d)] Assume $ x\ge 0$ is such that $y:= Ax \ne 0.$  Let $L := \text{span}(A_B)$ and decompose $y = y_L + y_{\perp}$
where $y_{\perp} = A^\perp_N x_N \in L^\perp$ and $y_L = A_B x_B +(A_N -  A^\perp_N) x_N \in L.$  Put $r:=\frac{\|y_\perp\|}{\|y\|} \in [0,1].$  Assume $r>0$ as otherwise $y = y_L\in \text{span}(A_B)$ and the statement holds with the better bound $\phi(A) \ge \vert\rho_B(A)\vert$ by proceeding exactly as in part (a).
Since $r >0$, we have  $x_N \ne 0$.  Put  $r_N:=\frac{\|y_\perp\|}{\|x_N\|_1}.$  From \eqref{rhoNdual} it follows that $r_N \ge \rho_N(A).$ 
Next, put $w:=\frac{1}{\|x_N\|_1}\left((A_N - A^{\perp}_N) x_N - y_L \right)$.  Observe that $\|w\| \le \dmax_{i\in N}\|a_i-a^\perp_i\| + \frac{\|y_L\|}{\|x_N\|_1} \le \|A\| + \frac{r_N\sqrt{1-r^2}}{r}$ and $w\in L$.  Hence by \eqref{rhoBprimal} there exists $\tilde x_B \ge 0, \|\tilde x_B\|_1 = 1$ such that
$A_B \tilde x_B = c w,$
where
\[
c:=\frac{\vert\rho_B(A)\vert r}{r\|A\| +  r_N\sqrt{1-r^2}} \in (0,1).
\]
Taking $\tilde x_N:= \frac{c}{\|x_N\|_1} x_N$ we get
\[
 A_N \tilde x_N - A_B \tilde x_B = \frac{c}{\|x_N\|_1} (y_\perp + y_L) = \frac{\vert\rho_B(A)\vert  r_N}{r\|A\| +  r_N\sqrt{1-r^2}} \frac{y}{\|y\|}.
\]
Thus letting $u := (1-c)x + (0,\tilde x_N),\; v = (\tilde x_B,0)$  we get $u,v\in \Delta_{n-1}, \; S(u) \subseteq S(x)$ and
\begin{equation}\label{case.c}
Au - Av = \left((1-c)\|Ax\| + \frac{\vert\rho_B(A)\vert  r_N}{r\|A\| +  r_N\sqrt{1-r^2}}\right) \frac{Ax}{\|Ax\|}.\end{equation}
Next, observe that
\begin{equation}\label{bound.case.c}
\begin{array}{rcl}
(1-c)\|Ax\| + \dfrac{\vert\rho_B(A)\vert  r_N}{r\|A\| +  r_N\sqrt{1-r^2}} 
&\ge & \dfrac{\vert\rho_B(A)\vert  r_N}{r\|A\| +  r_N\sqrt{1-r^2}} \\
&\ge &\dfrac{\vert\rho_B(A)\vert  r_N}{\sqrt{\|A\|^2 + r_N^2}} \\[3ex]&\ge &\dfrac{\vert\rho_B(A)\vert  \rho_N(A)}{\sqrt{\|A\|^2 + \rho_N(A)^2}}.
\end{array}
\end{equation}
The first inequality above follows because $c\in(0,1)$, the second one follows from 
\[
\dmax_{r\in[0,1]} \left(r\|A\| +  r_N\sqrt{1-r^2}\right) = \sqrt{\|A\|^2 +  r_N^2},\]
and the third one follows from $r_N \ge \rho_N(A)$.  
Putting~\eqref{case.c} and~\eqref{bound.case.c} together we get $\phi(A,x) \ge \dfrac{\vert\rho_B(A)\vert  \rho_N(A)}{\sqrt{\|A\|^2 + \rho_N(A)^2}}.$  
\qed
\end{description}

\begin{corollary}\label{bound.varrho} Assume $A = \matr{a_1 & \cdots & a_n} \in \R^{m\times n}$ is a nonzero matrix and $0\in \conv(A)$.  Then $w(A) \ge \phi(A) > 0$.
\end{corollary}
\proof  Apply Theorem~\ref{main.thm}.
Since $0 \in \conv(A)$, we have $B\ne \emptyset$ and thus case (b) cannot occur.  
If case (a) occurs then $\phi(A) \ge \vert\rho_B(A)\vert > 0$ since $\rho_B(A) <0$ as $L\ne \{0\}$.  If case (c) occurs then $\phi(A) \ge \rho_N(A) > 0$.  
Finally, if case (d) occurs then 
$\phi(A) \ge \dfrac{\vert\rho_B(A)\vert\rho_N(A)}{\sqrt{\|A\|^2 + \rho_N(A)^2}}>0,$ since both $\rho_N(A) > 0$ and $\rho_B(A) <0$ as $L\ne \{0\}$.  To finish, recall that $w(A) \ge \phi(A)$ as established in Section~\ref{sec.quadratic}.
\qed

\bigskip

We conclude with a few small examples that illustrate the values of $\phi(A),\vert\rho_B(A)\vert, \rho_N(A)$ and their connection with the bounds in Theorem~\ref{main.thm} for the three possible cases: $N = \emptyset, \; B = \emptyset,$ and  both $B,N \ne \emptyset.$

\begin{example}  Assume $\epsilon, \delta \in (0,1)$ and let $$
A = \matr{-1 & 1 & -1 & 1 & -1 & 1\\ 
-\epsilon & -\epsilon &\epsilon & \epsilon &\epsilon\delta & \epsilon\delta}.$$  In this case $B=\{1,2,3,4,5,6\}, \; N =\emptyset$.  It is easy to see that $\vert\rho_B(A)\vert = \epsilon$ and $\phi(A) = \phi(A,\bar x) = (1+\delta)\epsilon$ for $\bar x = \matr{0&  0 &0&0&1/2 & 1/2}\transp$.
\end{example}

\begin{example}  Assume $\delta \in (0,1)$ and let $A = \matr{1 & -1 \\ 
\delta & \delta}.$  In this case $B= \emptyset,\; N =\{1,2\}$.  It is easy to see that $\rho_N(A) = \delta$ and if we put $\bar A = \matr{A & 0}$ then $\phi(\bar A) = \phi(\bar A,\bar x) = \delta$ for $\bar x = \matr{1/2 & 1/2 & 0}\transp$.
\end{example}
   
\begin{example}  Assume $\epsilon,\delta \in (0,1)$ and let $$A = \matr{-1 & 1 & -1 & 1 & 0 & 0\\ 
-\epsilon & -\epsilon &\epsilon & \epsilon & 1 & -1 \\ 0 & 0 &0 & 0 & \delta & \delta}.$$  In this case $B=\{1,2,3,4\}, \; N =\{5,6\}$.  It is easy to see that $\vert\rho_B(A)\vert = \epsilon, \; \rho_N(A) = \delta.$ For $\bar x = \matr{0&  0 &\frac{1}{2(1+\epsilon)} & \frac{1}{2(1+\epsilon)} & 0 & \frac{\epsilon}{1+\epsilon}}\transp$ we get 
\[
A\bar x = \matr{0 \\ 0 \\ \dfrac{\epsilon \delta}{1+\epsilon}}.
\]
It thus follows that $\phi(A) \le \phi(A,\bar x) = \frac{2\epsilon \delta}{1+\epsilon}$.  On the other hand, Theorem~\ref{main.thm} implies that in this case 
$\phi(A) \ge \frac{\epsilon \delta}{\sqrt{\max(1+\epsilon^2,1+\delta^2)+\delta^2}}.$  In particular, $\epsilon\delta< \phi(A) < 2\epsilon\delta.$
\end{example}

\section*{Acknowledgements} 
We are grateful to Simon Lacoste-Julien and Martin Jaggi for their comments on a preliminary draft of this paper.  We are also grateful to two anonymous referees for their numerous constructive suggestions.  The first author's research has been supported by NSF grant CMMI-1534850.

\end{document}